# Sum-connectivity indices of trees and unicyclic graphs of fixed maximum degree

Zhibin Du[a], Bo Zhou[a]*, Nenad Trinajstić[b]

[a]Department of Mathematics, South China Normal University,

Guangzhou 510631, China

email: zhoubo@scnu.edu.cn

[b]The Rugjer Bošković Institute, P. O. Box 180, HR-10002 Zagreb, Croatia

email: trina@irb.hr

Abstract

We obtain the maximum sum-connectivity indices of graphs in the set of trees and in the set of unicyclic graphs respectively with given number of vertices and maximum degree, and determine the corresponding extremal graphs. Additionally, we deduce the $n$-vertex unicyclic graphs with the first two maximum sum-connectivity indices for $n \geq 4$.

*Corresponding author.

# 1. Introduction

The well-known Randić connectivity index, proposed by Randić [1] in 1975, is the most used molecular descriptor in quantitative structure property relationship (QSPR) and quantitative structure activity relationship [e.g., 1,4,5,6,11,14]. Mathematical properties of this index have also been studied [e.g., 2,7,9,10].

Let $G$ be a simple (molecular) graph with vertex set $V(G)$ and edge set $E(G)$. For $v \in V(G)$, $\Gamma(v)$ denotes the set of its (first) neighbors in $G$ and the degree of $v$ is $d_v = d_G(v) = |\Gamma(v)|$. The Randić connectivity index $R = R(G)$ of a graph $G$ is defined as [1]

$$R = R(G) = \sum_{uv \in E(G)} (d_u d_v)^{-1/2}.$$

A related connectivity index of a graph $G$, called the sum-connectivity index, is defined as [12]

$$\chi = \chi(G) = \sum_{uv \in E(G)} (d_u + d_v)^{-1/2}.$$

Correspondingly, we call $R(G)$ the product-connectivity index of $G$.

We determined in [15] the unique tree with fixed numbers of vertices and pendant vertices (vertices of degree one) with the minimum value of the sum-connectivity index, and trees with the minimum, second minimum and third minimum, and the maximum, second maximum and third maximum values of the sum-connectivity index, and we discussed its properties for a class of trees representing acyclic hydrocarbons.

The product- and sum-connectivity indices are highly intercorrelated quantities; for example, the value of the correlation coefficient is 0.99088 for 134 trees representing the lower alkanes from [3]. But, in general, the sum-connectivity index has a wider range [15].

Properties on the product-connectivity index for trees and unicyclic graphs may be found in [7,8,13]. Particularly, the product-connectivity index for trees with given maximum degree was investigated in [8].

For $2 \leq \Delta \leq n-1$, let $\mathbf{T}(n,\Delta)$ be the set of trees with $n$ vertices and maximum degree $\Delta$, and $\mathbf{U}(n,\Delta)$ the set of unicyclic graphs with $n$ vertices and maximum degree $\Delta$. Let $P_n$ and $C_n$ respectively be the path and the cycle on $n \geq 3$ vertices. In particular, $\mathbf{T}(n,2) = \{P_n\}$ and $\mathbf{U}(n,2) = \{C_n\}$.

In this paper, we obtain the maximum sum-connectivity indices of graphs in $\mathbf{T}(n,\Delta)$ and $\mathbf{U}(n,\Delta)$, respectively, and determine the corresponding extremal graphs. Recall that the $n$-vertex trees with the first two maximum sum-connectivity indices for $n \geq 4$ have been known [15]. From the result related to maximum degree for unicyclic graphs, we now deduce the $n$-vertex unicyclic graphs with the first two maximum sum-connectivity indices for $n \geq 4$.

## 2. Preliminaries

For a edge subset $E_1$ of the graph $G$ (the complement of $G$, respectively), $G - E_1$ ($G + E_1$, respectively) denotes the graph resulting from $G$ by deleting (adding, respectively) the edges in $E_1$.

**Lemma 1.** [15 For a connected graph $Q$ with at least two vertices and vertex $u \in V(Q)$, let $G_1$ be the graph obtained from $Q$ by attaching two paths $P_a$ and $P_b$ to $u$, and $G_2$ the graph obtained from $Q$ by attaching a path $P_{a+b}$ to $u$, where $a \geq b \geq 1$. Then $\chi(G_1) < \chi(G_2)$.

**Lemma 2.** For a connected graph $M$ with $|V(M)| \geq 3$, a vertex $u$ of degree two, let $H$ be the graph obtained from $M$ by attaching a path $P_a$ to $u$. Denote $u_1$ and $u_2$ by the two neighbors of $u$ in $M$, and $u'$ by the pendant vertex of the path attached to $u$ in $H$. If $\min\{d_H(u_1), d_H(u_2)\} \leq 4$, then for $H' = H - \{uu_2\} + \{u'u_2\}$, $\chi(H') > \chi(H)$.

**Proof.** Assume that $d_H(u_1) \geq d_H(u_2)$, and thus, $d_H(u_2) \leq 4$. If $d_H(u,u') = 1$, then

$$\chi(H') - \chi(H) = \left(\frac{1}{\sqrt{d_H(u_1)+2}} + \frac{1}{\sqrt{d_H(u_2)+2}} + \frac{1}{\sqrt{2+2}}\right) - \left(\frac{1}{\sqrt{d_H(u_1)+3}} + \frac{1}{\sqrt{d_H(u_2)+3}} + \frac{1}{\sqrt{1+3}}\right) > 0.$$

If $d_H(u,u') \geq 2$, then as $\frac{1}{\sqrt{x+2}} - \frac{1}{\sqrt{x+3}}$ is decreasing for $x > 0$, we have

$$\frac{1}{\sqrt{d_H(u_2)+2}} - \frac{1}{\sqrt{d_H(u_2)+3}} \geq \frac{1}{\sqrt{4+2}} - \frac{1}{\sqrt{4+3}} = \frac{1}{\sqrt{6}} - \frac{1}{\sqrt{7}}, \text{ and thus}$$

$$\begin{aligned}\chi(H') - \chi(H) &= \left( \frac{1}{\sqrt{d_H(u_1)+2}} + \frac{1}{\sqrt{d_H(u_2)+2}} + \frac{1}{\sqrt{2+2}} + \frac{1}{\sqrt{2+2}} \right) \\ &\quad - \left( \frac{1}{\sqrt{d_H(u_1)+3}} + \frac{1}{\sqrt{d_H(u_2)+3}} + \frac{1}{\sqrt{1+2}} + \frac{1}{\sqrt{2+3}} \right) \\ &= \left( \frac{1}{\sqrt{d_H(u_1)+2}} - \frac{1}{\sqrt{d_H(u_1)+3}} \right) + \left( \frac{1}{\sqrt{d_H(u_2)+2}} - \frac{1}{\sqrt{d_H(u_2)+3}} \right) + 1 - \frac{1}{\sqrt{3}} - \frac{1}{\sqrt{5}} \\ &> \left( \frac{1}{\sqrt{d_H(u_2)+2}} - \frac{1}{\sqrt{d_H(u_2)+3}} \right) + 1 - \frac{1}{\sqrt{3}} - \frac{1}{\sqrt{5}} \\ &\geq \left( \frac{1}{\sqrt{6}} - \frac{1}{\sqrt{7}} \right) + 1 - \frac{1}{\sqrt{3}} - \frac{1}{\sqrt{5}} > 0.\end{aligned}$$

The result follows. □

Attaching a path of length $r$ to a vertex $u$ of a graph $G$ means that adding an edge between $u$ and an terminal vertex of a path on $r$ vertices. Particularly, if $r = 1$, then a pendant vertex is attached. For $\frac{n}{2} \leq \Delta \leq n-1$, let $T_{n,\Delta}$ be the tree obtained by attaching $2\Delta+1-n$ pendant vertices and $n-\Delta-1$ paths of length two to a vertex. For $\frac{n+2}{2} \leq \Delta \leq n-1$, let $U_{n,\Delta}$ be the unicyclic graph obtained by attaching $2\Delta-n-1$ pendant vertices and $n-\Delta-1$ paths of length two to a vertex of a triangle.

## 3. Maximum Sum-Connectivity Index in $\mathbf{T}(n,\Delta)$

First we obtain the maximum sum-connectivity index of trees in $\mathbf{T}(n,\Delta)$ and determine the corresponding extremal graphs.

**Theorem 1.** Let $G \in \mathbf{T}(n,\Delta)$, where $2 \leq \Delta \leq n-1$. Then

$$\chi(G) \leq \begin{cases} \dfrac{2\Delta - n + 1}{\sqrt{\Delta+1}} + \dfrac{n-\Delta-1}{\sqrt{\Delta+2}} + \dfrac{n-\Delta-1}{\sqrt{3}} & \text{if } \dfrac{n}{2} \leq \Delta \leq n-1 \\ \dfrac{1}{2}(n-1-2\Delta) + \dfrac{\Delta}{\sqrt{3}} + \dfrac{\Delta}{\sqrt{\Delta+2}} & \text{if } 2 \leq \Delta \leq \dfrac{n-1}{2} \end{cases}$$

with equality if and only if $G = T_{n,\Delta}$ for $\dfrac{n}{2} \leq \Delta \leq n-1$, and $G$ is a tree obtained by attaching $\Delta$ paths of length at least two to a vertex for $2 \leq \Delta \leq \dfrac{n-1}{2}$.

**Proof.** The case $\Delta = 2$ is trivial. Suppose that $\Delta \geq 3$. Let $G$ be a tree in $\mathbf{T}(n,\Delta)$ with maximum sum-connectivity index. Let $v$ be a vertex of degree $\Delta$ in $G$. If there exists some vertex of degree more than two in $G$ different from $v$, then by Lemma 1, we may get a tree in $\mathbf{T}(n,\Delta)$ with larger sum-connectivity index, a contradiction. Thus, $v$ is the unique vertex of degree more than two in $G$. Let $k$ be the number of neighbors of $v$ with degree two. Then $k \leq \min\{n-\Delta-1, \Delta\}$. If $n-\Delta-1 \geq \Delta$, i.e., $\Delta \leq \dfrac{n-1}{2}$, then $1 \leq k \leq \Delta$. If $n-\Delta-1 < \Delta$, i.e., $\Delta \geq \dfrac{n}{2}$, then $0 \leq k \leq n-\Delta-1$. It is easily seen that

$$\dfrac{1}{\sqrt{3}} - \dfrac{1}{2} + \dfrac{1}{\sqrt{\Delta+2}} - \dfrac{1}{\sqrt{\Delta+1}} \geq \dfrac{1}{\sqrt{3}} - \dfrac{1}{2} + \dfrac{1}{\sqrt{3+2}} - \dfrac{1}{\sqrt{3+1}} > 0,$$

and then

$$\chi(G) = \frac{\Delta-k}{\sqrt{\Delta+1}} + \frac{k}{\sqrt{\Delta+2}} + \frac{k}{\sqrt{3}} + \frac{1}{2}(n-1-\Delta-k)$$

$$= \left(\frac{1}{\sqrt{3}} - \frac{1}{2} + \frac{1}{\sqrt{\Delta+2}} - \frac{1}{\sqrt{\Delta+1}}\right)k + \frac{1}{2}(n-1-\Delta) + \frac{\Delta}{\sqrt{\Delta+1}}$$

$$\leq \begin{cases} \left(\frac{1}{\sqrt{3}} - \frac{1}{2} + \frac{1}{\sqrt{\Delta+2}} - \frac{1}{\sqrt{\Delta+1}}\right)(n-\Delta-1) + \frac{1}{2}(n-1-\Delta) + \frac{\Delta}{\sqrt{\Delta+1}} & \text{if } \frac{n}{2} \leq \Delta \leq n-1 \\ \left(\frac{1}{\sqrt{3}} - \frac{1}{2} + \frac{1}{\sqrt{\Delta+2}} - \frac{1}{\sqrt{\Delta+1}}\right)\Delta + \frac{1}{2}(n-1-\Delta) + \frac{\Delta}{\sqrt{\Delta+1}} & \text{if } 2 \leq \Delta \leq \frac{n-1}{2} \end{cases}$$

$$= \begin{cases} \frac{2\Delta-n+1}{\sqrt{\Delta+1}} + \frac{n-\Delta-1}{\sqrt{\Delta+2}} + \frac{n-\Delta-1}{\sqrt{3}} & \text{if } \frac{n}{2} \leq \Delta \leq n-1 \\ \frac{1}{2}(n-1-2\Delta) + \frac{\Delta}{\sqrt{3}} + \frac{\Delta}{\sqrt{\Delta+2}} & \text{if } 2 \leq \Delta \leq \frac{n-1}{2} \end{cases}$$

with equality if and only if $k = n-\Delta-1$, i.e., each of the $n-\Delta-1$ neighbors of vertex $v$ of degree two is adjacent to a pendant vertex, i.e., $G = T_{n,\Delta}$ for $\frac{n}{2} \leq \Delta \leq n-1$, and $k = \Delta$, i.e., $G$ is a tree obtained by attaching $\Delta$ paths of length at least two to a vertex for $2 \leq \Delta \leq \frac{n-1}{2}$.

□

In Fig. 1, all the extremal graphs in Theorem 1 with $n = 7$ are given.

Fig. 1 comes here

**Fig. 1.** The 7-vertex trees with maximum sum-connectivity indices for $\Delta = 2, 3, 4, 5, 6$.

### 4. Maximum Sum-Connectivity Index in $\mathbf{U}(n,\Delta)$

Now we obtain the maximum sum-connectivity index of graphs in $\mathbf{U}(n,\Delta)$ and determine the corresponding extremal graphs. As a consequence, we deduce the *n*-vertex unicyclic graphs with the first and second maximum sum-connectivity indices for $n \geq 4$.

**Theorem 2.** Let $G \in \mathbf{U}(n,\Delta)$, where $2 \leq \Delta \leq n-1$. Then

$$\chi(G) \leq \begin{cases} \dfrac{n-\Delta-1}{\sqrt{3}} + \dfrac{n-\Delta+1}{\sqrt{\Delta+2}} + \dfrac{2\Delta-n-1}{\sqrt{\Delta+1}} + \dfrac{1}{2} & \text{if } \dfrac{n+2}{2} \leq \Delta \leq n-1 \\ \dfrac{\Delta-2}{\sqrt{3}} + \dfrac{\Delta}{\sqrt{\Delta+2}} + \dfrac{1}{2}(n-2\Delta+2) & \text{if } 2 \leq \Delta \leq \dfrac{n+1}{2} \end{cases}$$

with equality if and only if $G = U_{n,\Delta}$ for $\dfrac{n+2}{2} \leq \Delta \leq n-1$, and $G$ is a unicyclic graph obtained by attaching $\Delta-2$ paths of length at least two to a vertex of a cycle for $2 \leq \Delta \leq \dfrac{n+1}{2}$.

**Proof.** The case $\Delta = 2$ is trivial. Suppose that $\Delta \geq 3$, $G$ is a graph in $\mathbf{U}(n,\Delta)$ with maximum sum-connectivity index, and $C$ is the unique cycle of $G$. Let $v$ be a vertex of degree $\Delta$ in $G$.

First we consider $\Delta = 3$. If there is some vertex outside $C$ with degree three, then by Lemma 1, we may get a graph in $\mathbf{U}(n,3)$ with larger sum-connectivity index, a contradiction. If there are at least two vertices on $C$ with degree three, then by Lemma 2, we may get a graph in $\mathbf{U}(n,3)$ with larger sum-connectivity index, also a contradiction. Thus, $v \in V(C)$ and $v$ is the unique vertex in $G$ with degree three. Then either $\chi(G) = \dfrac{1}{2}(n-2) + \dfrac{2}{\sqrt{5}}$ when $v$ is adjacent to a vertex of degree one and two vertices of degree two for $n \geq 4$, or $\chi(G) = \dfrac{1}{2}(n-4) + \dfrac{3}{\sqrt{5}} + \dfrac{1}{\sqrt{3}}$ when $v$ is adjacent to three vertices of degree two for $n \geq 5$. Obviously, $\dfrac{1}{2}(n-2) + \dfrac{2}{\sqrt{5}} < \dfrac{1}{2}(n-4) + \dfrac{3}{\sqrt{5}} + \dfrac{1}{\sqrt{3}}$ for $n \geq 5$. Hence, $G$ is the graph obtained by attaching a pendant vertex to a triangle for $n = 4$, i.e., $G = U_{4,3}$, and a graph obtained by attaching a path of length at least two to a cycle for $n \geq 5$.

Now suppose that $\Delta \geq 4$. We will show that $v$ lies on $C$. Suppose that $v$ is not on $C$. Let $w$ be the vertex on $C$ such that $d_G(v,w) = \min\{d_G(v,x) : x \in V(C)\}$. If there is some vertex outside $C$ with degree more than two different from $v$, or if there is some vertex on $C$ with degree more than two different from $w$, then by Lemmas 1 and 2, we may get a graph in $\mathbf{U}(n,\Delta)$ with larger sum-connectivity index, a contradiction. Thus, $v$ and $w$ are the only vertices of degree more than two in $G$, and $d_G(v) = \Delta$, $d_G(w) = 3$ or $4$. Let $Q$ be the path connecting $v$ and $w$. Suppose that $d_G(w) = 4$. Denote $w_1$ by a neighbor of $w$ on $C$, $w_2$ by the neighbor of $w$ on $Q$, and $w'$ by the pendant vertex of the path attached to $w$. Consider $G_1 = G - \{ww_1\} + \{w'w_1\} \in U(n, \Delta)$. If $d_G(w, w') = 1$, then

$$\chi(G_1) - \chi(G) = \left( \frac{1}{\sqrt{d_G(w_2)+3}} + \frac{1}{\sqrt{2+3}} + \frac{1}{\sqrt{2+3}} + \frac{1}{\sqrt{2+2}} \right)$$
$$- \left( \frac{1}{\sqrt{d_G(w_2)+4}} + \frac{1}{\sqrt{1+4}} + \frac{1}{\sqrt{2+4}} + \frac{1}{\sqrt{2+4}} \right)$$
$$= \left( \frac{1}{\sqrt{d_G(w_2)+3}} - \frac{1}{\sqrt{d_G(w_2)+4}} \right) + \frac{1}{\sqrt{5}} + \frac{1}{2} - \frac{2}{\sqrt{6}} > 0.$$

If $d_G(w, w') \geq 2$, then

$$\chi(G_1) - \chi(G) = \left( \frac{1}{\sqrt{d_G(w_2)+3}} + \frac{1}{\sqrt{2+3}} + \frac{1}{\sqrt{2+3}} + \frac{1}{\sqrt{2+2}} + \frac{1}{\sqrt{2+2}} \right)$$
$$- \left( \frac{1}{\sqrt{d_G(w_2)+4}} + \frac{1}{\sqrt{2+4}} + \frac{1}{\sqrt{2+4}} + \frac{1}{\sqrt{1+2}} + \frac{1}{\sqrt{2+4}} \right)$$
$$= \left( \frac{1}{\sqrt{d_G(w_2)+3}} - \frac{1}{\sqrt{d_G(w_2)+4}} \right) + \frac{2}{\sqrt{5}} + 1 - \frac{1}{\sqrt{3}} - \frac{3}{\sqrt{6}} > 0.$$

In either case, $\chi(G_1) > \chi(G)$, a contradiction. Thus, $d_G(w) = 3$. Suppose that $v_1, v_2, \ldots, v_{\Delta-1}$ are the neighbors of $v$ outside $Q$. Let $d_i = d_G(v_i)$ for $i = 1, \ldots, \Delta - 1$. Note that $d_1, d_2 = 1$ or 2. Consider $G_2 = G - \{vv_3, \ldots, vv_{\Delta-1}\} + \{wv_3, \ldots, wv_{\Delta-1}\} \in U(n, \Delta)$. Note that $d_{G_2}(w) = \Delta$ and $d_{G_2}(v) = 3$. Then

$$\chi(G_2) - \chi(G) = \frac{1}{\sqrt{d_1+3}} - \frac{1}{\sqrt{d_1+\Delta}} + \frac{1}{\sqrt{d_2+3}} - \frac{1}{\sqrt{d_2+\Delta}} + \frac{2}{\sqrt{\Delta+2}} - \frac{2}{\sqrt{2+3}}$$
$$\geq \frac{1}{\sqrt{2+3}} - \frac{1}{\sqrt{2+\Delta}} + \frac{1}{\sqrt{2+3}} - \frac{1}{\sqrt{2+\Delta}} + \frac{2}{\sqrt{\Delta+2}} - \frac{2}{\sqrt{5}} = 0.$$

Since $d_{G_2}(v) = 3$, then by Lemma 1, we may get a graph $G'$ in $\mathbf{U}(n,\Delta)$ such that $\chi(G') > \chi(G_2) \geq \chi(G)$, a contradiction. Hence, we have shown that $v$ lies on $C$.

If there is some vertex outside $C$ with degree more than two, then by Lemma 1, we may get a graph in $\mathbf{U}(n,\Delta)$ with larger sum-connectivity index, a contradiction. Thus, $G$ is a graph obtained from $C$ by attaching $\Delta - 2$ paths to $v$, and attaching at most one path to some vertex on $C$ different from $v$. If there is some vertex on $C$ with degree three, then by Lemma 2, we may get a graph in $\mathbf{U}(n,\Delta)$ with larger sum-connectivity index, a contradiction. Thus, $G$ is a graph obtained from $C$ by attaching $\Delta - 2$ paths to $v$. Let $k$ be the number of neighbors of $v$ with degree two. Then $k \leq \min\{n-\Delta-1, \Delta-2\}$. If $n-\Delta-1 \geq \Delta-2$, i.e., $\Delta \leq \frac{n+1}{2}$, then $0 \leq k \leq \Delta - 2$. If $n-\Delta-1 < \Delta-2$, i.e., $\Delta \geq \frac{n+2}{2}$, then $0 \leq k \leq n-\Delta-1$. It is easily seen that

$$\frac{1}{\sqrt{3}} - \frac{1}{2} + \frac{1}{\sqrt{\Delta+2}} - \frac{1}{\sqrt{\Delta+1}} \geq \frac{1}{\sqrt{3}} - \frac{1}{2} + \frac{1}{\sqrt{4+2}} - \frac{1}{\sqrt{4+1}} > 0,$$

and then

$$\chi(G) = \frac{k}{\sqrt{3}} + \frac{k+2}{\sqrt{\Delta+2}} + \frac{\Delta-k-2}{\sqrt{\Delta+1}} + \frac{1}{2}(n-\Delta-k)$$

$$= \left(\frac{1}{\sqrt{3}} - \frac{1}{2} + \frac{1}{\sqrt{\Delta+2}} - \frac{1}{\sqrt{\Delta+1}}\right)k + \frac{\Delta-2}{\sqrt{\Delta+1}} + \frac{2}{\sqrt{\Delta+2}} + \frac{1}{2}(n-\Delta)$$

$$\leq \begin{cases} \left(\frac{1}{\sqrt{3}} - \frac{1}{2} + \frac{1}{\sqrt{\Delta+2}} - \frac{1}{\sqrt{\Delta+1}}\right)(n-\Delta-1) + \frac{\Delta-2}{\sqrt{\Delta+1}} + \frac{2}{\sqrt{\Delta+2}} + \frac{1}{2}(n-\Delta) \\ \qquad \text{if } \frac{n+2}{2} \leq \Delta \leq n-1 \\ \left(\frac{1}{\sqrt{3}} - \frac{1}{2} + \frac{1}{\sqrt{\Delta+2}} - \frac{1}{\sqrt{\Delta+1}}\right)(\Delta-2) + \frac{\Delta-2}{\sqrt{\Delta+1}} + \frac{2}{\sqrt{\Delta+2}} + \frac{1}{2}(n-\Delta) \\ \qquad \text{if } 2 \leq \Delta \leq \frac{n+1}{2} \end{cases}$$

$$= \begin{cases} \frac{n-\Delta-1}{\sqrt{3}} + \frac{n-\Delta+1}{\sqrt{\Delta+2}} + \frac{2\Delta-n-1}{\sqrt{\Delta+1}} + \frac{1}{2} & \text{if } \frac{n+2}{2} \leq \Delta \leq n-1 \\ \frac{\Delta-2}{\sqrt{3}} + \frac{\Delta}{\sqrt{\Delta+2}} + \frac{1}{2}(n-2\Delta+2) & \text{if } 2 \leq \Delta \leq \frac{n+1}{2} \end{cases}$$

with equality if and only if $k=n-\Delta-1$ for $\frac{n+2}{2}\leq\Delta\leq n-1$, i.e., $G=U_{n,\Delta}$, and $k=\Delta-2$ for $2\leq\Delta\leq\frac{n+1}{2}$, i.e., $G$ is a unicyclic graph obtained by attaching $\Delta-2$ paths of length at least two to a vertex of a cycle. □

In Fig. 2, all the extremal graphs in Theorem 2 with $n=7$ are given.

Fig. 2 comes here

**Fig. 2.** The 7-vertex unicyclic graphs with maximum sum-connectivity indices for $\Delta=2,3,4,5,6$ (for $\Delta=3$, there are three such graphs).

**Theorem 3.** Among the unicyclic graphs on $n\geq 4$ vertices, $C_n$ is the unique graph with maximum sum-connectivity index, which is equal to $\frac{n}{2}$, for $n=4$, $U_{4,3}$ is the unique graph with the second maximum sum-connectivity index, which is equal to $1+\frac{2}{\sqrt{5}}$, while for $n\geq 5$, the graphs obtained by attaching a path of length at least two to a cycle are the unique graphs with the second maximum sum-connectivity index, which is equal to $\frac{1}{2}(n-4)+\frac{1}{\sqrt{3}}+\frac{3}{\sqrt{5}}$.

**Proof.** The case $n=4$ may be checked directly. Suppose that $n\geq 5$ and $G$ is a unicyclic graph on $n$ vertices. Let $\Delta$ be the maximum degree of $G$, where $2\leq\Delta\leq n-1$. Let

$$f(x)=\frac{x-2}{\sqrt{3}}+\frac{x}{\sqrt{x+2}}+\frac{1}{2}(n-2x+2) \text{ for } x\geq 2.$$

If $\frac{n+2}{2}\leq\Delta\leq n-1$, then by Theorem 2,

$$\chi(G)\leq\frac{n-\Delta-1}{\sqrt{3}}+\frac{n-\Delta+1}{\sqrt{\Delta+2}}+\frac{2\Delta-n-1}{\sqrt{\Delta+1}}+\frac{1}{2}$$
$$=f(\Delta)+(n-2\Delta+1)\left(\frac{1}{\sqrt{3}}-\frac{1}{2}+\frac{1}{\sqrt{\Delta+2}}-\frac{1}{\sqrt{\Delta+1}}\right)<f(\Delta).$$

If $2 \leq \Delta \leq \frac{n+1}{2}$, then by Theorem 2, $\chi(G) \leq f(\Delta)$. Note that for $x \geq 2$, $\left(\frac{1}{2}x+2\right)(x+2)^{-3/2}$ is decreasing on $x \geq 2$, and then $\left(\frac{1}{2}x+2\right)(x+2)^{-3/2} \leq \left(\frac{1}{2}\cdot 2+2\right)(2+2)^{-3/2} = \frac{3}{8}$. Thus,

$$f'(x) = \frac{1}{\sqrt{3}} - 1 + \left(\frac{1}{2}x+2\right)(x+2)^{-3/2} \leq \frac{1}{\sqrt{3}} - 1 + \frac{3}{8} < 0,$$ and then $f(x)$ is decreasing for $x \geq 2$.

Now we have $\chi(G) < f(\Delta) < f(3) < f(2)$ for $3 < \frac{n+2}{2} \leq \Delta \leq n-1$, and $\chi(G) \leq f(\Delta) \leq f(3) < f(2)$ for $3 \leq \Delta \leq \frac{n+1}{2}$. Then obviously, $C_n$ is the unique $n$-vertex unicyclic graph with maximum sum-connectivity index $f(2)$, while the $n$-vertex unicyclic graphs with maximum degree three and sum-connectivity index $f(3)$ are the $n$-vertex unique graphs with the second maximum sum-connectivity index, and by Theorem 2, such graphs are the graphs obtained by attaching a path of length at least two to a cycle. □

**Acknowledgement** This work was supported by the National Natural Science Foundation of China (No. 11071089), the Guangdong Provincial Natural Science Foundation of China (No. S2011010005539) and the Ministry of Science, Education and Sports of Croatia (Grant No. 098-1770495-2919).